\input amstex
\documentstyle{amsppt}
%
\catcode`@=11
\redefine\output@{%
  \def\break{\penalty-\@M}\let\par\endgraf
  \ifodd\pageno\global\hoffset=105pt\else\global\hoffset=8pt\fi  
  \shipout\vbox{%
    \ifplain@
      \let\makeheadline\relax \let\makefootline\relax
    \else
      \iffirstpage@ \global\firstpage@false
        \let\rightheadline\frheadline
        \let\leftheadline\flheadline
      \else
        \ifrunheads@ 
        \else \let\makeheadline\relax
        \fi
      \fi
    \fi
    \makeheadline \pagebody \makefootline}%
  \advancepageno \ifnum\outputpenalty>-\@MM\else\dosupereject\fi
}
\catcode`\@=\active
\nopagenumbers
\def\negskp{\hskip -2pt}
\def\Ord{\operatorname{Ord}}
\def\Cl{\operatorname{Cl}}
\def\Rsim{\overset{\sssize R}\to\sim}
\def\blue#1{#1}

\catcode`#=11\def\diez{#}\catcode`#=6
\catcode`_=11\def\podcherkivanie{_}\catcode`_=8
\def\mycite#1{\cite{\blue{#1}}\immediate\special{ps:
     ShrHPSdict begin /ShrBORDERthickness 0 def}}

\def\mytag#1{%
    \tag#1}
\def\mythetag#1{\thetag{\blue{#1}}\immediate\special{ps:
     ShrHPSdict begin /ShrBORDERthickness 0 def}}
\def\myrefno#1{\no#1}
\def\myhref#1#2{\blue{#2}\immediate\special{ps:
     ShrHPSdict begin /ShrBORDERthickness 0 def}}

\def\mytheorem#1{\csname proclaim\endcsname{Theorem #1}}
\def\mytheoremwithtitle#1#2{\csname proclaim\endcsname{Theorem #1#2}}
\def\mythetheorem#1{\blue{#1}\immediate\special{ps:
     ShrHPSdict begin /ShrBORDERthickness 0 def}}
\def\mylemma#1{\csname proclaim\endcsname{Lemma #1}}
\def\mylemmawithtitle#1#2{\csname proclaim\endcsname{Lemma #1#2}}
\def\mythelemma#1{\blue{#1}\immediate\special{ps:
     ShrHPSdict begin /ShrBORDERthickness 0 def}}
\def\mycorollary#1{\csname proclaim\endcsname{Corollary #1}}

\def\mydefinition#1{\definition{Definition #1}}
\def\mythedefinition#1{\blue{#1}\immediate\special{ps:
     ShrHPSdict begin /ShrBORDERthickness 0 def}}

\pagewidth{360pt}
\pageheight{606pt}
\topmatter
\title
Transfinite normal and composition series of groups.
\endtitle
\author
R.~A.~Sharipov
\endauthor
\address 5 Rabochaya street, 450003 Ufa, Russia\newline
\vphantom{a}\kern 12pt Cell Phone: +7(917)476 93 48
\endaddress
\email \vtop to 30pt{\hsize=280pt\noindent
\myhref{mailto:r-sharipov\@mail.ru}
{r-sharipov\@mail.ru}\newline
\myhref{mailto:R\podcherkivanie Sharipov\@ic.bashedu.ru}
{R\_\hskip 1pt Sharipov\@ic.bashedu.ru}\vss}
\endemail
\urladdr
\vtop to 20pt{\hsize=280pt\noindent
\myhref{http://www.freetextbooks.narod.ru/index.html}
{http:/\negskp/www.freetextbooks.narod.ru/index.html}\newline
\myhref{http://sovlit2.narod.ru/}
{http:/\negskp/sovlit2.narod.ru}\vss}
\endurladdr
\abstract
    Normal and composition series of groups enumerated by ordinal numbers 
are studied. The Jordan-H\"older theorem for them is proved. 
\endabstract
\endtopmatter
\TagsOnRight
\document

\head
1. Introduction.
\endhead
     Finitely long normal and composition series of groups are well-known
in group theory. They are given by the following two definitions.
\mydefinition{1.1} Let $G$ be a group. A sequence of its subgroups 
$$
\hskip -2em
\{1\}=G_1\varsubsetneq G_2\varsubsetneq\ldots
\varsubsetneq G_n=G
\mytag{1.1}
$$
is called a normal series of subgroups for $G$ if the subgroup $G_i$ is a 
normal subgroup in $G_{i+1}$ for all $i=1,\,\ldots,\,n-1$.
\enddefinition
\mydefinition{1.2} A normal series \mythetag{1.1} is called a composition
series for $G$ if each factorgroup $G_{i+1}/G_i$ is simple.
\enddefinition
      According to the well-known Jordan-H\"older theorem (see \S\,3 of
Chapter~\uppercase\expandafter{\romannumeral 1} in \mycite{1}), 
the composition factors $G_{i+1}/G_i$ in a composition series of any group 
$G$ are unique up to some permutation of indices.
\mytheoremwithtitle{1.1}{ (Jordan-H\"older)} For any two composition series 
$$
\aligned
\{1\}=G_1\varsubsetneq G_2\varsubsetneq\ldots
\varsubsetneq G_n=G,\\
\vspace{1ex}
\{1\}=H_1\varsubsetneq H_2\varsubsetneq\ldots
\varsubsetneq H_m=G
\endaligned
$$
of a group $G$ their lengths are equal to each other, i\.\,e\. $n=m$, and 
there is some permutation $\sigma$ of the numbers $1,\,\ldots,\,n-1$ such 
that 
$$
G_{i+1}/G_i=H_{\sigma(i)+1}/H_{\sigma(i)}.
$$
\endproclaim
    In this paper we consider normal and composition series of the form 
\mythetag{1.1} where $n$ is some ordinal (transfinite) number. We prove 
the Jordan-H\"older theorem for such transfinite composition series. 
\head
2. Basic definitions.
\endhead
     Let's consider a transfinite series of subgroups of the form 
\mythetag{1.1} for some group $G$. Then $n$ is some ordinal number
(see Appendix~5 in \mycite{2}). Other indices in the series \mythetag{1.1} 
are ordinal numbers less than $n$. Let's recall that any ordinal number 
$\alpha$ is either a limit ordinal or a non-limit ordinal:
\roster
\rosteritemwd=0pt
\item"1)" if $\alpha$ is a limit ordinal, then for any $\beta<\alpha$ its 
successor $\beta+1$ is also less than $\alpha$;
\item"2)" if $\alpha$ is a non-limit ordinal, then $\alpha=\beta+1$ for some 
unique $\beta<\alpha$.
\endroster
\mydefinition{2.1} Let $G$ be a group. A transfinite sequence of subgroups 
$$
\hskip -2em
\{1\}=G_1\varsubsetneq G_2\varsubsetneq\ldots
\varsubsetneq G_n=G
\mytag{2.1}
$$
is called a transfinite normal series of subgroups for the group $G$ if the 
following two conditions are fulfilled:
\roster
\rosteritemwd=0pt
\item"1)" the subgroup $G_i$ is a 
normal subgroup in $G_{i+1}$ for any ordinal $i<n$;
\vskip 1ex
\item"2)" $G_\alpha=\dsize\bigcup_{\beta<\alpha}G_\beta$ for any limit ordinal 
$\alpha\leqslant n$.
\endroster
\enddefinition
\mydefinition{2.2} A group $G$ is called hypertranssimple if it has 
no normal series (nether finite nor transfinite) other than trivial one 
$\{1\}=G_1\varsubsetneq G_2=G$.
\enddefinition
\mydefinition{2.3} A transfinite normal series \mythetag{2.1} of a group $G$ 
is called a transfinite composition series of $G$ if for each ordinal $i<n$ 
the factorgroup $G_{i+1}/G_i$ is hypertranssimple.
\enddefinition
\head
3. The Jordan-H\"older theorem.
\endhead
\mydefinition{3.1} Let $G$ and $H$ be subgroups of some group. Then $G\cdot H$
is the subgroup generated by elements of both $G$ and $H$, i.\,e\. 
$G\cdot H=\langle G\cup H\rangle$.
\enddefinition
The subgroup $G\cdot H$ is composed by products of the form $(g_1\,h_1)\cdot
\ldots\cdot (g_n\,h_n)$ for some integer $n$, where $g_i\in G$ and $h_i\in H$
for all $i=1,\,\ldots,\,n$. If\/ $G$ or $H$ is a normal subgroup in a group
enclosing both $G$ and $H$, then 
$$
(g_1\,h_1)\cdot\ldots\cdot (g_n\,h_n)=g\,h
$$
for some $g\in G$ and some $h\in H$. In this case $G\cdot H$ is composed by
products $g\,h$, where $g\in G$ and $h\in H$.\par
     Note that $G\cdot H=H\cdot G$ even if $G$ and $H$ are not subgroups of
an Abelian group. Indeed, $G\cup H=H\cup G$. For this reason $G\cdot H
=\langle G\cup H\rangle=\langle H\cup G\rangle=H\cdot G$.
\mylemmawithtitle{3.1}{ (Zassenhaus)} Let\/ $\tilde G$ and $\tilde H$ be 
subgroups of some group and let $G$ and $H$ be normal subgroups 
of\/ $\tilde G$ and $\tilde H$ respectively. Then $G\cdot (\tilde G\cap H)$ 
is a normal subgroup in $G\cdot (\tilde G\cap\tilde H)$ and 
$H\cdot (\tilde H\cap G)$ is a normal subgroup in 
$H\cdot (\tilde H\cap\tilde G)$. The corresponding factorgroups are 
isomorphic, i\.\,e\.
$$
(G\cdot (\tilde G\cap\tilde H))/(G\cdot (\tilde G\cap H))\cong 
(H\cdot (\tilde H\cap\tilde G))/(H\cdot (\tilde H\cap G)).
$$
\endproclaim
The lemma~\mythelemma{3.1} is also known as the butterfly lemma. Its proof 
can be found in \S\,3 of Chapter~\uppercase\expandafter{\romannumeral 1} 
in \mycite{1}.
\mydefinition{3.2} A transfinite normal series $\{1\}=\tilde G_1\varsubsetneq 
\tilde G_2\varsubsetneq\ldots\varsubsetneq\tilde G_p=G$ is called a refinement
for a transfinite normal series $\{1\}=G_1\varsubsetneq G_2\varsubsetneq\ldots
\varsubsetneq G_n=G$ if each subgroup $G_i$ coincides with some subgroup 
$\tilde G_j$.
\enddefinition
\mydefinition{3.3} Two transfinite normal series $\{1\}=G_1\varsubsetneq G_2
\varsubsetneq\ldots\varsubsetneq G_n=G$ and $\{1\}=H_1\varsubsetneq H_2
\varsubsetneq\ldots\varsubsetneq H_m=G$ of a group $G$ are called isomorphic 
if there is a one-to-one mapping that associates each ordinal number $i<n$
with some ordinal number $j<m$ in such a way that 
$G_{i+1}/G_i\cong H_{j+1}/H_j$.
\enddefinition
\mytheorem{3.1} Arbitrary two transfinite normal series of a group\/ $G$ have 
isomorphic refinements. 
\endproclaim
     Let $\{1\}=G_1\varsubsetneq G_2\varsubsetneq\ldots\varsubsetneq 
G_n=G$ and $\{1\}=H_1\varsubsetneq H_2\varsubsetneq\ldots\varsubsetneq H_m=G$
be two transfinite normal series of some group $G$. The subgroups $G_i$ in 
the first series are indexed by ordinal numbers $i\leqslant n$, the subgroups 
$H_j$ in the  second series are indexed by ordinal numbers $j\leqslant m$. 
Let's consider the following subgroups of $G$ indexed by two indices $i$ and 
$j$ being ordinal numbers:
$$
\hskip -2em
\aligned
&G_{ij}=G_i\cdot(G_{i+1}\cap H_j)\text{, \ where \ }i<n\text{\ \ and \ }
j\leqslant m\,\\
\vspace{1ex}
&H_{\kern -0.4pt j\kern 0.4pt i}=H_j\kern 0.4pt\cdot(H_{j+1}\cap G_i)
\text{, \ where \ }i\leqslant n\text{\ \ and \ }j<m.
\endaligned
\mytag{3.1}
$$
Applying the Zassenhaus butterfly lemma~\mythelemma{3.1} to the subgroups
\mythetag{3.1}, we find that $G_{ij}$ is a normal subgroup in $G_{i\,j+1}$,
$H_{\kern -0.4pt j\kern 0.4pt i}$ is a normal subgroup in 
$H_{\kern -0.4pt j\kern 0.4pt\,i+1}$, and
$$
\hskip -2em
G_{i\,j+1}/G_{ij}\cong H_{\kern -0.4pt j\kern 0.4pt\,i+1}/H_{\kern -0.4pt j
\kern 0.4pt i}.
\mytag{3.2}
$$
The isomorphism \mythetag{3.2} is a base for proving the 
theorem~\mythetheorem{3.1}. However, it is not a proof since the subgroups 
$G_{ij}$ and $H_{\kern -0.4pt j\kern 0.4pt\,i}$ do not form transfinite 
normal series yet. In order to complete our proof we need some auxiliary 
lemmas.
\mydefinition{3.4} An ordered set $I$ \,is called totally ordered or 
linearly ordered if any two elements $i_1$ and $i_2$ of $I$ are comparable, 
i\.\,e\. $i_1\neq i_2$ means $i_1<i_2$ or $i_2<i_1$. 
\enddefinition
\mydefinition{3.5} A linearly ordered set $I$ is called well ordered if 
every non-empty subset $A\subseteq I$ has a minimal element $a_{\min}\in A$. 
It is clear that such a minimal element $a_{\min}$ in $A$ is unique.
\enddefinition
\mylemma{3.2} Let $I$ and $J$ be two well ordered sets. If we denote through 
$I\times J$ the set of pairs $(i,j)$, where $i\in I$ and $j\in J$, and if we 
equip $I\times J$ with the lexicographic ordering, then $I\times J$ is also 
a well ordered set.
\endproclaim
\demo{Proof} The lexicographic ordering in $I\times J$ means that 
$(i_1,j_1)<(i_2,j_2)$ if $i_1<i_2$ or if $i_1=i_2$ and $j_1<j_2$. 
It is easy to see that any two pairs $(i_1,j_1)$ and $(i_2,j_2)$ are 
comparable in this lexicographic ordering. Indeed, if $i_1\neq i_2$, 
since $I$ is linearly ordered, we have $i_1<i_2$ or $i_2<i_1$. The
inequality $i_1<i_2$ implies $(i_1,j_1)<(i_2,j_2)$, the converse 
inequality $i_2<i_1$ implies $(i_2,j_2)<(i_1,j_1)$. In both of these 
cases the pairs $(i_1,j_1)$ and $(i_2,j_2)$ are comparable.\par
      If $i_1=i_2$ and $(i_1,j_1)\neq (i_2,j_2)$, then, since $J$ is linearly 
ordered, we have $j_1<j_2$ or $j_2<j_1$. In the first case the equality 
$i_1=i_2$ and the inequality $j_1<j_2$ lead to $(i_1,j_1)<(i_2,j_2)$. 
Otherwise, from $i_2=i_1$ and from $j_2<j_1$ we derive 
$(i_2,j_2)<(i_1,j_1)$. Again, the pairs $(i_1,j_1)$ and $(i_2,j_2)$ 
appear to be comparable.\par
      Thus, the set $I\times J$ with the lexicographic ordering is linearly 
ordered. Let $A\subseteq I\times J$ be some non-empty subset of $I\times J$. 
Let's denote through $A_I$ the projection of $A$ onto the first component 
of the direct product $I\times J$:
$$
\hskip -2em
A_I=\{i\in I\!:\,\exists\,j\in J\text{\ such that }(i,j)\in A\}.
$$
The subset $A_I\subset I$ is not empty. Since $I$\ \,is well ordered, there 
is a unique minimal element $i_{\min}\in A_I$. Now let's denote trough $A_J$ 
the following subset of $J$:
$$
A_J=\{j\in J\!:\,(i_{\min},j)\in A\}.
$$
The subset $A_J\subseteq J$ is also not empty. Since $J$ is well ordered, 
the subset $A_J$ has a unique minimal element $j_{\min}\in A_J$. Then 
$(i_{\min},j_{\min})$ is a unique minimal element of the subset 
$A\subseteq I\times J$. The lemma~\mythelemma{3.2} is proved.
\qed\enddemo
\mydefinition{3.6} Let $\Rsim$ be an equivalence relation in a linearly 
ordered set $I$. The equivalence relation $\Rsim$ is said to be concordant 
with the linear ordering in $I$\/ if\/ $i_1<i_2<i_3$ and $i_1\Rsim i_3$ 
imply $i_1\Rsim i_2$ and $i_2\Rsim i_3$.
\enddefinition
     Let $I$ be a linearly ordered set and let $\Rsim$ be an equivalence 
relation in $I$ concordant with its linear ordering. Then the factorset 
$I/R$ can be equipped with the factorordering. It is introduced as follows: 
for any two distinct equivalence classes $\Cl_R(i_1)\neq\Cl_R(i_2)$ we 
declare $\Cl_R(i_1)<\Cl_R(i_2)$ if $i_1<i_2$. This definition is 
self-consistent, i\.\,e\. the inequality $\Cl_R(i_1)<\Cl_R(i_2)$ does not 
depend on our choice of the representatives $i_i$ and $i_2$ within the 
equivalence classes $\Cl_R(i_1)$ and $\Cl_R(i_2)$. Indeed, assume that 
$\Cl_R(i_3)=\Cl_R(i_1)$ and $i_3\neq i_1$. Since $I$\/ is linearly ordered, 
we have two options: $i_3<i_2$ or $i_2<i_3$. If $i_2<i_3$ we would have 
$i_1<i_2<i_3$ and $i_1\Rsim i_3$. Applying the concordance condition (see 
Definition~\mythedefinition{3.6}), we would derive $i_1\Rsim i_2$, which 
contradicts $\Cl_R(i_1)\neq\Cl_R(i_2)$. Thus, the option $i_2<i_3$ is 
excluded and we have $i_3<i_2$, which means $\Cl_R(i_3)<\Cl_R(i_2)$.\par
      Similarly, if we assume that $\Cl_R(i_4)=\Cl_R(i_2)$ and $i_4\neq i_2$, 
then, applying the concordance condition, we derive $i_1<i_4$, which means 
$\Cl_R(i_1)<\Cl_R(i_4)$. Thus, the definition of the factorordering in the 
factorset $I/R$ is self-consistent.
\mylemma{3.3} If\/ $I$ is a well ordered set and if\/ $\Rsim$ is an 
equivalence relation concordant with the ordering in $I$, then the factorset 
$I/R$ equipped with the factorordering is a well ordered set too.
\endproclaim
\demo{Proof} The well ordered set $I$ is linearly ordered. It is easy to see 
that the factorordering in $I/R$ is a linear ordering too. Let $A$ be a 
non-empty subset of the factorset $I/R$. Let's denote through $\tilde A$ 
its preimage in $I$ 
$$
\tilde A=\{i\in I\!:\,\Cl_R(i)\in A\}.
$$
The subset $\tilde A\subset I$ is non-empty. Hence, since $I$ is well ordered,
there is a unique minimal element $\tilde a_{\min}$ in $\tilde A$. Let's 
denote $a=\Cl_R(\tilde a_{\min})$. It is clear that $a$ is a minimal element 
of $A$. The proof is over.\qed
\enddemo
     Let's return back to the subgroups \mythetag{3.1}. In order to describe 
these subgroups let's introduce the following sets of ordinal numbers:
$$
\xalignat 2
&I=\{\alpha\in\Ord\!:\,\alpha\leqslant n\},
&&I'=\{\alpha\in\Ord\!:\,\alpha<n\},\\
&J=\{\beta\in\Ord\!:\,\beta\leqslant m\},
&&J'=\{\beta\in\Ord\!:\,\beta<m\}.
\endxalignat
$$
The subgroups $G_{ij}$ in \mythetag{3.1} are indexed by the set $I'\times J$ 
ordered lexicographically. Similarly, the subgroups $H_{\kern -0.4pt j
\kern 0.4pt i}$ are indexed by the set $J'\times I$ ordered lexicographically.
\par
     There are two types of elements in the ordered set $I'\times J$ --- 
regular elements of the form $(i,j)$, where $i<n$ and $j<m$, and irregular 
elements of the form $(i,m)$, where $i<n$. Similarly, the ordered set 
$J'\times I$ has regular elements of the form $(\kern 0.1pt j\kern 0.1pt,
\kern 0.5pt i)$, where $i<n$ and $j<m$, and irregular elements 
$(\kern 0.1pt j\kern 0.1pt ,\kern 0.5pt n)$, where $j<m$. Note that regular 
elements of the set $I'\times J$ are in a one-to-one association with regular 
elements of the set $J'\times I$. Indeed, we have 
$$
\hskip -2em
\CD
@>\theta>>\\
\vspace{-4ex}
(i,j)@.
(\kern 0.1pt j\kern 0.1pt,\kern 0.5pt i).\\
\vspace{-4.2ex}
@<<\theta< 
\endCD
\mytag{3.3}
$$
The formula \mythetag{3.3} means that the association $\theta$ is composed 
by two bijective mappings inverse to each other. Both of them are denoted 
through $\theta$.
\par
     Let $\alpha=(i,j)$ be a regular element of the well ordered set 
$I'\times J$ and let $\beta=(\kern 0.1pt j\kern 0.1pt,\kern 0.5pt i)$ be 
its associated element in $J'\times I$. Then 
$$
\xalignat 2
&\beta=\theta(\alpha),
&&\alpha=\theta(\beta).
\endxalignat
$$
The next element to $\alpha$ is $\alpha+1=(i,j+1)$. Similarly, the next 
element to $\beta$ is $\beta+1=(\kern 0.1pt j\kern 0.1pt,\kern 0.5pt i+1)$. 
Applying the formula \mythetag{3.2}, we derive
$$
G_{\alpha+1}/G_\alpha\cong H_{\beta+1}/H_\beta
$$
for each regular element $\alpha\in I'\times J$ and its associated regular 
element $\beta\in J'\times I$.\par
     Now let $\alpha=(i,m)$ be an irregular element of the well ordered set 
$I'\times J$. Then, applying the first formula \mythetag{3.1}, we get 
$G_\alpha=G_{i\kern 0.5pt m}=G_{i+1}$. The next element to $\alpha$ in this 
case is $\alpha+1=(i+1,1)$. Applying the first formula \mythetag{3.1} again, 
we get $G_{\alpha+1}=G_{i+1\kern 1pt 1}=G_{i+1}$. Thus we have
$$
\hskip -2em
G_{\alpha+1}=G_\alpha
\mytag{3.4}
$$
for each irregular element $\alpha\in I'\times J$. Similarly, the second 
formula \mythetag{3.1} yields
$$
\hskip -2em
H_{\beta+1}=H_\beta
\mytag{3.5}
$$
for each irregular element $\beta\in J'\times I$.\par
     If $n$ is a non-limit ordinal, then $n=(n-1)+1$. In this case the 
well ordered set $I'\times J$ has the maximal element $\alpha_{\max}=(n-1,m)$ 
which is irregular. If $n$ is a limit ordinal, then the well ordered set 
$I'\times J$ has no maximal elements. In this case we extend this set by 
adjoining an auxiliary element $\alpha_{\max}=(n,1)$ to it:
$$
\hskip -2em
IJ=\cases\kern 2.5em I'\times J &\text{if $n$ is a non-limit ordinal;}\\
(I'\times J)\cup\{(n,1)\}&\text{if $n$ is a limit ordinal.}
\endcases
\mytag{3.6}
$$
We extend the ordering of $I'\times J$ to $IJ$ by declaring the auxiliary 
element $(n,1)$ to be greater than all elements of $I'\times J$. It is easy 
to see that the extended set \mythetag{3.6} is a well ordered set possessing 
the maximal element $\alpha_{\max}$. If $n$ is a non-limit ordinal, the first 
formula \mythetag{3.1} yields
$$
\hskip -2em
G_{\alpha_{\max}}=G_{n-1\kern 1pt m}=G.
\mytag{3.7}
$$
If $n$ is a limit ordinal, we extend the first formula \mythetag{3.1} by 
setting
$$
\hskip -2em
G_{\alpha_{\max}}=G_{n\kern 1pt 1}=G.
\mytag{3.8}
$$\par
     The second well ordered set $J'\times I$ can also be incomplete. In 
this case we extend it by adjoining an auxiliary element to this set:
$$
\hskip -2em
JI=\cases\kern 2.5em J'\times I &\text{if $m$ is a non-limit ordinal;}\\
(J'\times I)\cup\{(m,1)\}&\text{if $m$ is a limit ordinal.}
\endcases
\mytag{3.9}
$$
The ordering of the extended set \mythetag{3.9} extends the ordering of 
$J'\times I$ so that the auxiliary element $\beta_{\max}=(m,1)$ is the 
maximal element of $JI$. We extend the second formula \mythetag{3.1} to 
this element by setting
$$
\hskip -2em
H_{\beta_{\max}}=H_{m\kern 1pt 1}=G.
\mytag{3.10}
$$
If $m$ is a non-limit ordinal the maximal element of $JI$ is given by the 
formula
$$
\hskip -2em
H_{\beta_{\max}}=H_{m-1\kern 1pt n}=G.
\mytag{3.11}
$$
Like the extended set $IJ$ in \mythetag{3.6}, the extended set $JI$ in 
\mythetag{3.9} is a well ordered set possessing its maximal element 
$\beta_{\max}$.\par
     Due to \mythetag{3.4} and \mythetag{3.5} irregular elements of the 
well ordered sets $IJ$ and $JI$ can be identified (glued) with the regular 
elements next to them. As for the regular elements, some of them can also 
be identified with other regular elements. In order to perform such 
identifications we introduce some special equivalence relations in the 
well ordered sets $IJ$ and $JI$. 
\mydefinition{3.7} Two elements $\alpha_1=(i_1,j_1)$ and $\alpha_2=(i_2,j_2)$ 
of the set $IJ$ are declared to be equivalent $\alpha_1\Rsim\alpha_2$ if
$G_{\alpha_1}=G_{\alpha_2}$, i\.\,e\. if $G_{i_1j_1}=G_{i_2j_2}$. 
\enddefinition
\mydefinition{3.8} Two elements $\beta_1=(\kern 0.1pt j_1,\kern 0.5pt i_1)$ 
and $\beta_2=(\kern 0.1pt j_2,\kern 0.5pt i_2)$ of the set $JI$ are declared 
to be equivalent $\beta_1\Rsim\beta_2$ if $G_{\beta_1}=G_{\beta_2}$, 
i\.\,e\. if $G_{j_1i_1}=G_{j_2i_2}$. 
\enddefinition
\mylemma{3.4} The equivalence relations introduced by the 
definitions~\mythedefinition{3.7} and \mythedefinition{3.8} are concordant
with the orderings in $IJ$ and $JI$ in the sense of the 
definition~\mythedefinition{3.6}.
\endproclaim
\demo{Proof} Using the formulas \mythetag{3.1} and its extensions 
\mythetag{3.8} and \mythetag{3.10}, it is easy to prove that for two 
elements $\alpha_1$ and $\alpha_2$ of the ordered set $IJ$\/ the 
inequality $\alpha_1<\alpha_2$ implies 
$G_{\alpha_1}\subseteq G_{\alpha_2}$. Relying on the 
definition~\mythedefinition{3.6}, assume that $\alpha_1<\alpha_2<\alpha_3$ 
and $\alpha_1\Rsim\alpha_3$. Then we have the following relationships:
$$
\xalignat 2
&\hskip -2em
G_{\alpha_1}\subseteq G_{\alpha_2}\subseteq G_{\alpha_3}
&&G_{\alpha_1}=G_{\alpha_3}.\quad
\mytag{3.12}
\endxalignat
$$
From \mythetag{3.12} we immediately derive $G_{\alpha_1}=G_{\alpha_2}$ and
$G_{\alpha_2}=G_{\alpha_3}$ which means that $\alpha_1\Rsim\alpha_2$ and
$\alpha_2\Rsim\alpha_3$. Thus the proposition of the lemma~\mythelemma{3.4} 
is proved for the equivalence relation $\Rsim$ in $IJ$. In the case of the 
second ordered set $JI$ the proof is similar since in this case the 
inequality $\beta_1<\beta_2$ implies $H_{\beta_1}\subseteq H_{\beta_2}$.
\qed\enddemo
     Factorizing the well ordered sets $IJ$ and $JI$ with respect to the 
equivalence relations introduced in the definitions~\mythedefinition{3.7}
and \mythedefinition{3.8}, we get two factorsets
$$
\xalignat 2
&\hskip -2em
IJR=(IJ)/R,
&&JIR=(JI)/R.
\mytag{3.13}
\endxalignat
$$
According to the lemmas~\mythelemma{3.3} and \mythelemma{3.4}, the factorsets 
\mythetag{3.13} are well ordered. It is known that each well ordered set is 
isomorphic to some initial segment in the class of ordinal numbers (see 
Appendix~3 in \mycite{2}). The well ordered sets $IJ$ and $JI$ have the 
maximal elements $\alpha_{\max}$ and $\beta_{\max}$. Their equivalence classes 
are maximal elements in the factorsets $IJR$ and $JIR$ respectively. For this 
reason there are two ordinal numbers $p$ and $q$ such that
$$
\hskip -2em
\aligned
&IJR\cong \{1,\,\ldots,\,p\kern 0.5pt\}
=\{\alpha\in\Ord\!:\,\alpha\leqslant p\kern 0.5pt\},\\
\vspace{1ex}
&JIR\cong \{1,\,\ldots,\,q\kern 0.5pt\}
=\{\beta\in\Ord\!:\,\beta\leqslant q\kern 0.5pt\}.
\endaligned
\mytag{3.14}
$$
The ordinal numbers $p$ and $q$ in \mythetag{3.14} correspond to the 
equivalence classes of the maximal elements $\alpha_{\max}$ and 
$\beta_{\max}$ respectively. 
\par
      According to the definitions~\mythedefinition{3.7} and 
\mythedefinition{3.8}, the factorsets \mythetag{3.13} can be used for 
indexing the subgroups $G_{ij}$ and $H_{\kern -0.4pt j\kern 0.4pt i}$. 
Applying \mythetag{3.14} and taking into account \mythetag{3.7}, 
\mythetag{3.8}, \mythetag{3.10}, and \mythetag{3.11}, we get two 
transfinite sequences of subgroups of $G$:
$$
\hskip -2em
\aligned
&{1}=G_1\varsubsetneq\ldots\varsubsetneq G_p=G,\\
\vspace{1ex}
&{1}=H_1\varsubsetneq\ldots\varsubsetneq H_q=G.
\endaligned
\mytag{3.15}
$$
\par
     The next step is to prove that the series \mythetag{3.15} are 
transfinite normal series in the sense of the definition~\mythedefinition{2.1}. 
Let's consider the first series \mythetag{3.15}. Assume that $r$ be some 
ordinal less than $p$. Then $r+1\leqslant p$ is the ordinal next to $r$. The 
ordinals $r$ and $r+1$ are associated with two immediately adjacent classes 
$A$ and $A+1$ 
$$
\xalignat 2
&\hskip -2em
A=\Cl_R(a), 
&&A+1=\Cl_R(b)
\mytag{3.16}
\endxalignat
$$
in the factorset $IJR$. The upper class $A+1$ in \mythetag{3.16} is a 
non-empty subset of the well ordered set $IJ$ (as well as the lower class 
$A$). It has some unique minimal element $b_{\kern 0.1pt\min}$. Without 
loss of generality we can assume that $b=b_{\kern 0.1pt\min}$.
\mylemma{3.5} For any two immediately adjacent classes $A$ and $A+1$ of\/ 
the factorset $IJR$ the upper class $A+1$ has the minimal element 
$b_{\kern 0.1pt\min}$, while the lower class $A$ has the maximal element 
$a_{\kern 0.1pt\max}$ such that $b_{\kern 0.1pt\min}=a_{\kern 0.1pt\max}+1$. 
\endproclaim
\demo{Proof} In order to prove the lemma~\mythelemma{3.5} it is sufficient to
prove that $b_{\kern 0.1pt\min}$ is a non-limit element of the well ordered 
set $IJ$. The proof is by contradiction. Assume that $b_{\kern 0.1pt\min}$ is
a limit element of the set $IJ$. There are two options for 
$b_{\kern 0.1pt\min}$:
\roster
\rosteritemwd=5pt
\item"1)" $b_{\kern 0.1pt\min}=(i,1)$, where $i\leqslant n$ is a limit ordinal;
\item"2)" $b_{\kern 0.1pt\min}=(i,j)$, where $i<n$ is an arbitrary ordinal and 
$j\leqslant m$ is a limit ordinal.
\endroster
In the first case, applying the first formula \mythetag{3.1} or the formula 
\mythetag{3.8}, we derive 
$$
\hskip -2em
G_{b_{\kern 0.1pt\min}}=G_{i\kern 1pt 1}=G_i.
\mytag{3.17}
$$
Since $i$ is a limit ordinal, applying the item 2 of the definition
\mythedefinition{2.1}, we get
$$
\hskip -2em
G_i=\bigcup_{\alpha<i}G_{\alpha}=\bigcup_{\alpha<i}G_{\alpha+1}
=\bigcup_{\alpha<i}G_{\alpha\kern 1pt m}.
\mytag{3.18}
$$
Combining the formula \mythetag{3.17} with the formula \mythetag{3.18}, 
we derive 
$$
\hskip -2em
G_{b_{\kern 0.1pt\min}}=\bigcup_{\alpha<i}G_{\alpha\kern 1pt m}.
\mytag{3.19}
$$
Note that $\alpha<i$ implies the inequality $(\alpha,m)<(i,1)$, where $(i,1)
=b_{\kern 0.1pt\min}$. Since $b_{\kern 0.1pt\min}$ is the minimal element of 
its class, we have $\Cl_R(\alpha,m)<\Cl_R(b_{\kern 0.1pt\min})=A+1$. This
inequality can be rewritten as follows:
$$
\hskip -2em
\Cl_R(\alpha,m)\leqslant A\text{\ \ for all \ }\alpha<i.
\mytag{3.20}
$$
The class $A$ is not empty. Therefore there is at least one ordinal number
$\alpha<i$ such that $\Cl_R(\alpha,m)=A$. This equality and the above
inequality \mythetag{3.20} lead to the following relationships for subgroups:
$$
\hskip -2em
\aligned
&G_{\alpha\kern 1pt m}\subseteq G_A\text{\ \ for all \ }\alpha<i;\\
\vspace{1ex}
&G_{\alpha\kern 1pt m}=G_A\text{\ \ for some \ }\alpha<i.
\endaligned
\mytag{3.21}
$$
Applying \mythetag{3.21} to \mythetag{3.19}, we derive 
$G_{b_{\kern 0.1pt\min}}\!=G_A$. But, on the other hand, we have 
$G_{b_{\kern 0.1pt\min}}\!=G_{A+1}\neq G_A$ since $A\neq A+1$. 
Thus, the first of the above two options for $b_{\kern 0.1pt\min}$ leads 
to a contradiction.\par
     Let's proceed to the second option. In this case $b_{\kern 0.1pt\min}
=(i,j)$, where $i<n$ is some arbitrary ordinal and $j\leqslant m$ is some 
limit ordinal. The item 2 of the definition~\mythedefinition{2.1} applied
to the series $\{1\}=H_1\varsubsetneq\ldots\varsubsetneq H_m=G$ says that 
$$
\hskip -2em
H_j=\bigcup_{\beta<j}H_\beta.
\mytag{3.22}
$$
Combining \mythetag{3.22} with the first formula \mythetag{3.1}, we derive
the formula 
$$
\hskip -2em
G_{b_{\kern 0.1pt\min}}=G_{ij}=\bigcup_{\beta<j}G_{i\beta}.
\mytag{3.23}
$$
This formula \mythetag{3.23} is similar to \mythetag{3.19}. Using the
arguments quite similar to the above ones, we get the following relationships:
$$
\hskip -2em
\aligned
&G_{i\beta}\subseteq G_A\text{\ \ for all \ }\beta<j;\\
\vspace{1ex}
&G_{i\beta}=G_A\text{\ \ for some \ }\beta<j.
\endaligned
\mytag{3.24}
$$
Applying \mythetag{3.24} to \mythetag{3.23}, we derive 
$G_{b_{\kern 0.1pt\min}}\!=G_A$. But, on the other hand, we have 
$G_{b_{\kern 0.1pt\min}}\!=G_{A+1}\neq G_A$. As we see, the second 
of the above two options for $b_{\kern 0.1pt\min}$ also leads to a 
contradiction. Due to these two contradictions we conclude that 
$b_{\kern 0.1pt\min}$ is a non-limit element of the well ordered
set $IJ$. Then $b_{\kern 0.1pt\min}=a+1$ for some unique element 
$a\in IJ$. It is clear that $a=a_{\kern 0.1pt\max}$ is the maximal 
element in the lower class $A$. The lemma is proved.
\qed\enddemo
\mylemma{3.6} For any two immediately adjacent classes $B$ and $B+1$ of\/ 
the factorset $JIR$ the upper class $B+1$ has the minimal element 
$b_{\kern 0.1pt\min}$, while the lower class $B$ has the maximal element 
$a_{\kern 0.1pt\max}$ such that $b_{\kern 0.1pt\min}=a_{\kern 0.1pt\max}+1$. 
\endproclaim
     The lemma~\mythelemma{3.6} is an analog of the previous 
lemma~\mythelemma{3.5}.For this reason it does not require a separate 
proof.\par
{\bf A remark}. Not each subset in a well ordered set has a maximal 
element, but, according to the lemma~\mythelemma{3.5}, each equivalence 
class in $IJ$\/ has. According the lemma~\mythelemma{3.6}, the same is 
true for equivalence classes in $JI$.\par
\mylemma{3.7} Let $A$ be a limit class of the factorset $IJR$. Then
$$
\hskip -2em
G_A=\bigcup_{B<A}G_B.
\mytag{3.25}
$$ 
\endproclaim
\demo{Proof} Let $a_{\kern 0.1pt\min}$ be the minimal element of the class 
$A$. Then $a_{\kern 0.1pt\min}$ is a limit element of the well ordered set 
$IJ$. There are the following options for $a_{\kern 0.1pt\min}$:
\roster
\rosteritemwd=5pt
\item"1)" $a_{\kern 0.1pt\min}=(i,1)$, where $i\leqslant n$ is a limit ordinal;
\item"2)" $a_{\kern 0.1pt\min}=(i,j)$, where $i<n$ is an arbitrary ordinal and 
$j\leqslant m$ is a limit ordinal.
\endroster
In the fist case we have the following formula for the subgroup $G_A$:
$$
\hskip -2em
G_A=G_{a_{\kern 0.1pt\min}}=\bigcup_{\alpha<i}G_{\alpha}.
\mytag{3.26}
$$
The arguments are the same as in deriving the formulas \mythetag{3.17} and
\mythetag{3.18}. Note that each class $B<A$ in this case is represented by 
some element $b=(\alpha,j)$, where $\alpha<i$ and $1\leqslant j\leqslant m$. 
The inequalities $1\leqslant j\leqslant m$ yield
$$
\hskip -2em
G_\alpha=G_{\alpha 1}\subseteq G_{\alpha j}\subseteq G_{\alpha\kern 1pt m}
=G_{\alpha+1}.
\mytag{3.27}
$$
Since $i$ is a limit ordinal, combining \mythetag{3.26} and \mythetag{3.27}, 
we derive 
$$
G_A=\bigcup\Sb \alpha<i\\j\leqslant m\endSb G_{\alpha j}
=\bigcup_{B<A}G_B.
$$
Thus, the formula \mythetag{3.25} is proved for the first case where 
$a_{\kern 0.1pt\min}=(i,1)$ and $i\leqslant n$ is some limit ordinal.\par
     Let's proceed to the second case $a_{\kern 0.1pt\min}=(i,j)$, where $i<n$ 
is some arbitrary ordinal and where $j\leqslant m$ is some limit ordinal. 
In this case we have
$$
\hskip -2em
G_{a_{\kern 0.1pt\min}}=G_{ij}=\bigcup_{\beta<j}G_{i\kern 0.5pt\beta}.
\mytag{3.28}
$$
The formula \mythetag{3.28} is identical to \mythetag{3.23}. Each class $B<A$ 
is represented by some element $b=(\alpha,\beta)$ such that 
$b<a_{\kern 0.1pt\min}$. Due to the lexicographic ordering in $IJ$ there 
are two types of such elements $b=(\alpha,\beta)$ --- those with $\alpha<i$ 
and those with $\alpha=i$ and $\beta<j$. For the elements of the first type 
we have
$$
\hskip -2em
G_{\alpha\beta}\subseteq G_{\alpha\kern 1pt m}=G_{\alpha+1\kern 1pt 1}
\subseteq G_{i\kern 1pt 1}.
\mytag{3.29}
$$
Due to the inclusions \mythetag{3.29} we can derive the following formula:
$$
\hskip -2em
\bigcup_{B<A}G_B=\biggl(\,\bigcup\Sb \alpha<i\\\beta\leqslant m\endSb 
G_{\alpha\beta}\biggr)\cup\Bigl(\bigcup_{\beta<j}G_{i\kern 0.5pt\beta}
\Bigr)=\bigcup_{\beta<j}G_{i\kern 0.5pt\beta}.
\mytag{3.30}
$$
Comparing \mythetag{3.30} with \mythetag{3.28}, we see that the formula 
\mythetag{3.25} is proved for the second case. Thus the lemma~\mythelemma{3.7}
is completely proved.
\qed\enddemo
\mylemma{3.8} Let $B$ be a limit class of the factorset $JIR$. Then
$$
\hskip -2em
G_B=\bigcup_{A<B}G_A.
\mytag{3.31}
$$ 
\endproclaim
The lemma~\mythelemma{3.8} is an analog of the lemma~\mythelemma{3.7}, while
the formula \mythetag{3.31} is an analog of the formula \mythetag{3.25}.
\par
     Now let's return back to the series \mythetag{3.15}. Recall that they 
were built by the subgroups \mythetag{3.1}, \mythetag{3.8} and \mythetag{3.10} 
in such a way that they obey the item 1 of the 
definition~\mythedefinition{2.1}. On the other hand, the 
lemmas~\mythelemma{3.7} and \mythelemma{3.8} prove that they obey the item 
2 of the definition~\mythedefinition{2.1} as well. Thus the series 
\mythetag{3.15} are transfinite normal series of subgroups for the group 
$G$ in the sense of the definition~\mythedefinition{2.1}.\par
     Note that $G_{i\kern 1pt 1}=G_i$ and $H_{\kern -0.4pt j\kern 1pt 1}
=H_{\kern -0.4pt j}$. For this reason all subgroups of the initial series 
$\{1\}=G_1\varsubsetneq\ldots\varsubsetneq G_n=G$ and $\{1\}=H_1\varsubsetneq
\ldots\varsubsetneq H_m=G$ are among the subgroups of the series 
\mythetag{3.15}, i.\,e\. the transfinite normal series \mythetag{3.15} are 
refinements of the initial series $\{1\}=G_1\varsubsetneq\ldots
\varsubsetneq G_n=G$ and $\{1\}=H_1\varsubsetneq\ldots\varsubsetneq H_m=G$ 
in the sense of the definition~\mythedefinition{3.2}.
\mylemma{3.9} The transfinite normal series \mythetag{3.15} of the group $G$ 
are isomorphic to each other in the sense of the 
definition~\mythedefinition{3.3}.
\endproclaim
\demo{Proof} The proof is based on the lemmas~\mythelemma{3.5} and 
\mythelemma{3.6}. Due to \mythetag{3.14} each factorgroup $G_{A+1}/G_A$ of 
the first normal series \mythetag{3.15} is associated with some pair of 
immediately adjacent equivalence classes $A$ and $A+1$ being the elements 
of the factorset $IJR=(IJ)/R$. According to the lemma~\mythelemma{3.5}, 
the lower of these two classes has the maximal element $a_{\kern 0.1pt\max}
=(i,j)$. Then 
$$
\hskip -2em
G_{i\,j+1}/G_{ij}=G_{A+1}/G_A\neq\{1\}.
\mytag{3.32}
$$
Applying the formula \mythetag{3.2} to \mythetag{3.32}, we get
$$
\hskip -2em
G_{A+1}/G_A\cong H_{\kern -0.4pt j\kern 0.4pt\,i+1}/H_{\kern -0.4pt j
\kern 0.4pt i}\neq\{1\}.
\mytag{3.33}
$$
Let's denote through $B$ the equivalence class of the element 
$b=(\kern 0.1pt j\kern 0.1pt,\kern 0.5pt i)$, i\.\,e\. let 
$B=\Cl_R(\kern 0.1pt j\kern 0.1pt,\kern 0.5pt i)$. Note that $i<n$ and 
$j<m$ in \mythetag{3.33}. Hence $b+1=(\kern 0.1pt j\kern 0.1pt,\kern 0.5pt 
i+1)$ and
$$
\hskip -2em
H_{b+1}/H_b=H_{\kern -0.4pt j\kern 0.4pt\,i+1}/H_{\kern -0.4pt j
\kern 0.4pt i}\neq\{1\}.
\mytag{3.34}
$$
The formula \mythetag{3.34} means that $b=b_{\kern 0.1pt\max}$ is the maximal
element of its class $B$ and 
$$
\hskip -2em
G_{A+1}/G_A\cong H_{B+1}/H_B\neq\{1\}.
\mytag{3.35}
$$\par
     Thus, we have proved that the upper mapping $\theta$ in \mythetag{3.3} 
goes through the factorization procedures in $IJ$ and $JI$ and associates each 
class $A$ having its successor $A+1$ in $IJR$ with some unique class $B$
having its successor $B+1$ in $JIR$. Applying the above arguments to a class
$B$ of $JIR$, we find that same is true for the lower mapping $\theta$ in 
\mythetag{3.3}. Like the initial mappings \mythetag{3.3}, the factorized 
mappings $\theta$ are bijective and inverse to each other. The formula 
\mythetag{3.35} shows that these factorized mappings constitute an isomorphism 
for the transfinite normal series \mythetag{3.15}. The lemma~\mythelemma{3.9} 
\pagebreak is proved.
\qed\enddemo
     The lemma~\mythelemma{3.9} proved just above completes the proof of
the theorem~\mythetheorem{3.1}. 
\mylemma{3.10} If\/ $\{1\}=G_1\varsubsetneq\ldots \varsubsetneq G_n=G$ 
is a transfinite composition series of a group $G$, then it has no refinements 
different from itself. 
\endproclaim
\demo{Proof} The proof is by contradiction. Assume that the transfinite 
composition series $\{1\}=G_1\varsubsetneq\ldots\varsubsetneq G_n=G$ of the
group $G$ has some nontrivial refinement $\{1\}=\tilde G_1\varsubsetneq\ldots
\varsubsetneq\tilde G_p=G$ different from $\{1\}=G_1\varsubsetneq\ldots
\varsubsetneq G_n=G$. According to the definition~\mythedefinition{3.3},
each $G_i$ of the first series coincides with some $G_j$ in the second series.
In other words we have an injective mapping
$$
\sigma\!:\,\{i\in\Ord\!:\, i\leqslant n\}\longrightarrow
\{j\in\Ord\!:\, j\leqslant p\}
$$
such that $\tilde G_{\sigma(i)}=G_i$ for each $i\leqslant n$. Since 
$G_1=\tilde G_1$ and $G_n=\tilde G_p,$, we have 
$$
\xalignat 2
&\hskip -2em
\sigma(1)=1,
&&\sigma(n)=p\neq n.
\mytag{3.36}
\endxalignat
$$
The set $I=\{i\in\Ord\!:\, i\leqslant n\}$ subdivides into two subsets 
$$
\xalignat 2
&\hskip -2em
I_1=\{i\in I\!:\,\sigma(i)=i\},
&&I_2=\{i\in I\!:\,\sigma(i)\neq i\}.
\quad
\mytag{3.37}
\endxalignat
$$
Due to \mythetag{3.36} both subsets \mythetag{3.37} are not empty. Let's
denote through $i_2=i_{\min}$ the minimal element of the subset $I_2$.
\par
     The index $i_2$ is a non-limit ordinal. Indeed, otherwise we would have
$$
\xalignat 2
&G_{i_2}=\bigcup_{\alpha<i_2}G_\alpha,
&&\tilde G_{i_2}=\bigcup_{\alpha<i_2}\tilde G_\alpha
=\bigcup_{\alpha<i_2}G_\alpha,
\endxalignat
$$
which yields $\tilde G_{i_2}=G_{i_2}$ and $\sigma(i_2)=i_2$. Since $i_2$ is a
non-limit ordinal, the first subset $I_1$ in \mythetag{3.37} has the maximal
element $i_1=i_{\max}$ such that $i_2=i_1+1$. In other words, $i_1$ and $i_2$ 
are two neighboring ordinals such that 
$$
\xalignat 2
&\hskip -2em
\sigma(i_1)=i_1,
&&\sigma(i_2)>i_2.
\mytag{3.38}
\endxalignat
$$
Due to the relationships \mythetag{3.38} we have the following segment in the 
refined normal series $\{1\}=\tilde G_1\varsubsetneq\ldots\varsubsetneq
\tilde G_p=G$ of the group $G$:
$$
\hskip -2em
G_{i_1}=\tilde G_{i_1}\varsubsetneq\tilde G_{i_1+1}\varsubsetneq
\ldots\varsubsetneq\tilde G_{i_1+q}=G_{i_2}
\mytag{3.39}
$$
Here $q>1$ is some ordinal number. Note that $i_2=i_1+1$. Hence $G_{i_1}$ 
is a normal subgroup of $G_{i_2}$. Then $G_{i_1}$ is a normal subgroup of 
each group in the sequence \mythetag{3.39}. Passing to factorgroups in 
the sequence \mythetag{3.39}, we get
$$
\hskip -2em
\{1\}=\tilde G_{i_1}/G_{i_1}\varsubsetneq
\ldots\varsubsetneq\tilde G_{i_1+q}/G_{i_1}=G_{i_1+1}/G_{i_1}.
\mytag{3.40}
$$
It is easy to show that \mythetag{3.40} is a normal series of the factorgroup
$G_{i_1+1}/G_{i_1}$. But, since $\{1\}=G_1\varsubsetneq\ldots\varsubsetneq 
G_n=G$ is a transfinite composition series, its factorgroup 
$G_{i_1+1}/G_{i_1}$ is hypertranssimple (see Definition~\mythedefinition{2.3}). 
Applying the definition~\mythedefinition{2.2}, we find that the normal series 
\mythetag{3.40} should be trivial, i\.\,e\. it should look like 
$$
\pagebreak
\{1\}=\tilde G_{i_1}/G_{i_1}\varsubsetneq\tilde G_{i_1+1}/G_{i_1}
=G_{i_1+1}/G_{i_1}.
$$
This means that $q=1$ and $\tilde G_{i_1+1}=G_{i_1+1}$. Then 
$\tilde G_{i_2}=G_{i_2}$ and $\sigma(i_2)=i_2$ which contradicts 
\mythetag{3.38}. The contradiction obtained completes the proof.
\qed\enddemo
\mytheoremwithtitle{3.2}{ (Jordan-H\"older)}Any two transfinite composition 
series of a group $G$ are isomorphic. 
\endproclaim
     Proving the theorem~\mythetheorem{3.2} is the main goal of this paper. 
Now its proof is immediate from the theorem~\mythetheorem{3.1} and the
lemma~\mythelemma{3.10}.\par
{\bf A remark}. Note that if we have two transfinite composition series
$$
\aligned
\{1\}=G_1\varsubsetneq G_2\varsubsetneq\ldots
\varsubsetneq G_n=G,\\
\vspace{1ex}
\{1\}=H_1\varsubsetneq H_2\varsubsetneq\ldots
\varsubsetneq H_m=G,
\endaligned
$$
their isomorphism (see Definition~\mythedefinition{3.3}) does not mean
$n=m$. It means only that $n$ and $m$ are two ordinal numbers of the
same cardinality, i\.\,e\. $|n|=|m|$. 
\head
4. Dedicatory.
\endhead
This paper is dedicated to my aunt Halila Yusupova who was born before
the First World War and lived her long life in works and permanent cares.

\enddocument
\end